\documentclass[12pt]{amsart}
\usepackage{amsmath,amsthm,amsfonts,amssymb,times,latexsym,mathabx,url}

\newtheorem{theorem}{Theorem}[section]

\newtheorem{lem}[theorem]{Lemma}

\newtheorem{rem}[theorem]{Remark}
\numberwithin{equation}{section}

\renewcommand{\L}{\mathcal{L}}
\newcommand{\e}{\varepsilon}

\renewcommand{\a}{\alpha}
\renewcommand{\b}{\beta}

\newcommand{\N}{\mathbb{N}}
\newcommand{\Z}{\mathbb{Z}}
\renewcommand{\leq}{\leqslant}
\renewcommand{\geq}{\geqslant}

\renewcommand{\L}{\Lambda}

\makeatletter
\renewcommand{\pmod}[1]{\allowbreak\mkern7mu({\operator@font mod}\,\,#1)}
\makeatother

\newcommand{\C}{\mathbb C}

\newcommand{\T}{\mathbb T}
\newcommand{\R}{\mathbb R}
\renewcommand{\L}{\Lambda}

\renewcommand{\a}{\alpha}
\renewcommand{\b}{\beta}
\newcommand{\g}{\gamma}
\renewcommand{\leq}{\leqslant}
\renewcommand{\geq}{\geqslant}

\begin{document}

\title[Trigonometric polynomials with frequencies in the set of cubes]{Trigonometric polynomials with frequencies \\ in the set of cubes}
\author{Mikhail R. Gabdullin, Sergei V. Konyagin}
\date{}

\address{Department of mathematics, 1409 West Green Street, University of Illinois at Urbana-Champaign, Urbana, IL 61801, USA; Steklov Mathematical Institute,
	Gubkina str., 8, Moscow, 119991, Russia}
\email{gabdullin.mikhail@yandex.ru, mikhailg@illinois.edu}

\address{Lomonosov Moscow State University, Leninskie Gory str., 1, Moscow, Russia, 119991; Steklov Mathematical Institute,
	Gubkina str., 8, Moscow, 119991, Russia}
\email{konyagin23@mi-ras.ru}



\begin{abstract}
We prove that for any $\e>0$ and any trigonometric polynomial $f$ with frequencies in the set $\{n^3: N \leq n\leq N+N^{2/3-\e}\}$, one has
$$
\|f\|_4 \ll \e^{-1/4}\|f\|_2 
$$
with implied constant being absolute. We also show that the set $\{n^3: N\leq n\leq N+(0.5N)^{1/2}\}$ is a Sidon set.
\end{abstract}

\date{\today}

\maketitle

\section{Introduction} 

We say that a trigonometric polynomial $f$ has frequencies in a set $A\subseteq\Z$ if, for some $a_n\in\C$, $f(x)=\sum_{n\in A}a_ne(nx)$ (here and in what follows $e(x)=e^{2\pi ix}$).  Let $p>2$. Recall that a set $A\subseteq \Z$ is said to be a $\L_p$-set if there exists a constant $C(A,p)>0$ such that for any trigonometric polynomial $f$ with frequencies in the set $A$, the inequality
\begin{equation}\label{1.1}
\|f\|_p \leq C(A,p)\|f\|_2
\end{equation}
holds, where $\|f\|_p=\left(\int_{\T}|f(x)|^pdx\right)^{1/p}$ and $\T=\R/\Z$. It is easy to show that if a set $A$ is a $\L_p$-set for some $p>2$, then the counting function of $A$ cannot grow too fast. To be more precise, let $A_N=A\cap[-N,N]$ and $f(x)=\sum_{n\in A_N}e(nx)$; then by the properties of the Dirichlet kernel $D_N)=\sum_{|k|\leq N}e(kx)$, H\"older's inequality and (\ref{1.1}),
$$
|A_N|=f(0)=(f*D_N)(0)=\int_{\T}f(-y)D_N(y)dy \leq \|f\|_p\|D_N\|_q \ll |A_N|^{1/2}N^{1/p}
$$
(here $q$ is defined by $1/p+1/q=1$), and, hence, $|A_N| \ll N^{2/p}$. Using probabilistic methods, Bourgain \cite{Bour} showed that this upper bound is tight: for instance, Theorem 3 in \cite{Bour} states that for any $p>2$, there exists a $\L_p$-subset $A$ of the set of prime numbers of maximal density, that is, with $\varliminf_{N\to\infty}|A_N|N^{-2/p}\gg1$. However, that proof does not produce any examples of $\L_p$-sets. On the other hand, we note that there is the classical result that any lacunary sequence is a $\L_p$-set for every $p>2$ (see, for example, \cite{Gr}, Theorem 3.6.4). 

In this regard, it is of a great interest to find polynomial sequences which are $\L_p$-sets. There is a famous and still unsolved conjecture that the set of squares $\{n^2: n\in \N\}$ is a $\L_p$-set for any $2<p<4$, which was discussed by W.~Rudin (see the end of section 4.6 in \cite{Rud}).
On the other hand, it is well-known (see the inequality (1.7) of \cite{Bour}) that
\begin{equation}\label{1.2}
\left\|\sum_{n\leq N}e(n^2x)\right\|_4 \asymp N^{1/2}(\log N)^{1/4},
\end{equation}
so the set of squares is not a $\L_4$-set. However, one can consider the sequence of squares in shorter intervals: a conjecture of Cilleruelo and C\'ordoba \cite{CC} asserts that, for any $\g\in(0,1)$ and any trigonometric polynomial $f$ with frequencies in the set $\{n^2: N\leq n\leq N+N^{\g}\}$, the inequality 
\begin{equation}\label{1.3}
\|f\|_4 \leq c(\g)\|f\|_2	
\end{equation}
holds with some $c(\g)>0$ depending only on $\g$. Recently the first author has made a progress on this conjecture, proving it for any $\g<(\sqrt5-1)/2=0.618\ldots$; see \cite{G}. We refer the interested reader to the works \cite{CG} and \cite{G} for a more detailed discussion of the case of squares.

In this paper we focus our attention on trigonometric polynomials with frequencies in the set of cubes. Since this set is even more sparse than that of squares, it is natural to expect same $\L_p$-properties for it as well. Moreover, there is no obstruction of the type (\ref{1.2}) for cubes due to the fact that $\|\sum_{n=1}^Ne(n^3x)\|_4 \asymp N^{1/2}$; this inequality follows immediately from the result of C.~Hooley \cite{Hoo} (it is also shown in the proof of Theorem 2.6 of \cite{Nat} that the number of nontrivial solutions of $a_1^3+b_1^3=a_2^3+b_2^3$ in $0\leq a_i,b_i\leq N$ is $O(N^{5/3+\e})$ for any $\e>0$; see also \cite{Woo} for an overview of bounds for cubic exponential sums). Thus, unlike the case of squares, it is reasonable to conjecture that the set of cubes $\{n^3: n\in\N \}$ is a $\L_4$-set. However, a little is known: even showing that
$$
\left\|\sum_{n=1}^Na_ne(n^3x)\right\|_4 \ll \left(\sum_{n=1}^N|a_n|^2\right)^{1/2}\!\!\!\cdot(\log N)^{O(1)}
$$
seems to be beyond the reach of current technique. Nevertheless, if one considers cubes in short intervals as in Cilleruelo-Cord\'oba conjecture (\ref{1.3}), then we can say something stronger than the corresponding known bounds for squares. Our main result is the following.  

\begin{theorem}\label{th1.1}
For any $\e>0$ and any trigonometric polynomial $f$ with frequencies in the set $\{n^3: N \leq n\leq N+N^{2/3-\e}\}$, 
$$
\|f\|_4 \ll \e^{-1/4}\|f\|_2. 
$$
\end{theorem}

We emphasize that $\L_p$-property is intimately connected to the arithmetic properties of a set $A\subset \Z$. For example, recall that a set $A\subset\Z$ is called a Sidon set if any $m\in\Z$ has at most one representation in the form $m=a+b$ with $a,b\in A$. It is easy to show that any Sidon set is a $\L_4$-set (see Lemma \ref{lem2.1} below) and that Sidon property is stronger in general. The set of cubes is not a Sidon set: there is a famous example of Ramanujan $1^3+12^3=9^3+10^3$, and, moreover, it is well-known that the number of representations of $m$ in the form $m=a^3+b^3$, where $a$ and $b$ are positive integers, is an unbounded function of $m$ (see Proposition 5.3 in \cite{Sil}).   

To compare $\L_4$ and Sidon properties for our case of cubes in a short interval, we also provide the following estimate (which is an analog of Proposition 4.1 in \cite{G}). 

\begin{theorem}\label{th1.2}
The set $\{n^3: N\leq n\leq N+(0.5N)^{1/2}\}$ is a Sidon set. Moreover, this statement is sharp up to the constant $(0.5)^{1/2}$.
\end{theorem}

Finally, we remark that it is natural to expect that the sets $\{n^k: n\in\N\}$ are $\L_4$-sets for all $k\geq4$ as well. Note also that there is an open conjecture of P.~Erd\H{o}s which asserts that the set of fifth powers $\{n^5: n\in\N \}$ is a Sidon set (again, probably the same holds for higher powers). We do not follow these questions in this paper.

\bigskip 

\textbf{Notation.} We use Vinogradov's $\ll$ notation: both $F\ll G$ and $F=O(G)$ mean that there exists a constant $C>0$ such that $|F|\leq CG$. We write $F\asymp G$ if $G\ll F\ll G$.

\bigskip 

\textbf{Acknowledgements.} This research was carried out at Lomonosov Moscow State University with the financial support of the Russian Science Foundation (grant no. 22-11-00129). 

\section{Proof of Theorem \ref{th1.1}}\label{sec2} 

We will rely on the following simple estimate, which is the inequality (6.1) of \cite{CG} (see also Lemma 2.1 in \cite{G} and the following Remark 2.3).  

\begin{lem}\label{lem2.1}
Let $A$ be a finite set of integers and
$$
r_A^+(m)=\#\{(n_1,n_2)\in A\times A: n_1+n_2=m\}.
$$
Then for $f(x)=\sum_{n\in A}a_n e(nx)$, we have
\begin{equation}\label{2.1}
\|f\|_4 \leq \left(\max_{m\in\Z} r_A^+(m)\right)^{1/4}\!\|f\|_2 \, .
\end{equation}
\end{lem}

Now let $N$ be large enough and $k=N^{2/3-\e}$. Assume that a number $m$ is represented as $m=u^3+v^3$ with $N\leq u,v\leq N+k$. Then
\begin{equation}\label{2.2}
m=u^3+v^3=(u+v)\Big((u+v)^2+3(u-v)^2\Big)/4,
\end{equation}
and, hence, $0\leq 4m-(u+v)^3 \leq 3(u-v)^2(u+v) \leq 6(N+k)k^2$. Further, since $u+v\geq 2N$ and $m\geq 2N^3$,  
$$
0\leq (4m)^{1/3}-(u+v)=\frac{4m-(u+v)^3}{(4m)^{2/3}+(4m)^{1/3}(u+v)+(u+v)^2} \leq \frac{6(N+k)k^2}{12N^2}<k^2N^{-1}. 
$$
Note that, for a fixed $m$, the value of $u+v$ determines the representation of the form (\ref{2.2}). Thus, denoting $A=\{n^3: N\leq n \leq N+k\}$ and taking into account that $u+v$ is a divisor of $m$, we see that 
\begin{multline}\label{2.3}
r_A^+(m)\leq \#\big\{d|4m: (4m)^{1/3}-k^2N^{-1} \leq d\leq (4m)^{1/3}\big\}\leq \\ \#\left\{d|4m: (4m)^{2/3} \leq d\leq\frac{4m}{(4m)^{1/3}-k^2N^{-1}}\right\}.	
\end{multline}
Since $k^2N^{-1}=N^{1/3-2\e}\leq m^{1/9-2\e/3}$ and $m$ is large, we obtain
\begin{multline}\label{2.4}
\frac{4m}{(4m)^{1/3}-k^2N^{-1}} \leq (4m)^{2/3}(1-4^{-1/3}m^{-2/9-2\e/3})^{-1} \leq \\ (4m)^{2/3}(1+m^{-2/9-2\e/3})\leq (4m)^{2/3}+(4m)^{4/9-2\e/3}. 
\end{multline}
Now we need the following result, which is actually \cite[Corollary 3.8]{CC} and can also be proven by a different method from \cite[Theorem 3.1]{G}. 

\begin{theorem}\label{th2.2}
Let $0<\a<1$ and $0<\b<\a^2$. For any positive integer $m$, we have 
$$
\#\big\{d|m: m^{\a}\leq d\leq m^{\a}+m^{\b}\big\}\ll (\a^2-\b)^{-1}
$$
with the implied constant being absolute.
\end{theorem}

Theorem \ref{2.2}, together with (\ref{2.3}) and (\ref{2.4}), implies that 
$$
r_A^+(m)\leq \#\big\{d|4m: (4m)^{2/3}\leq d\leq (4m)^{2/3}+(4m)^{4/9-2\e/3}\big\}\ll \e^{-1},
$$
and Theorem \ref{th1.1} follows from Lemma \ref{lem2.1}.

\begin{rem}
It is natural to ask whether the conclusion of Theorem \ref{th2.2} holds for $0<\b<\a$. To be more precise, we make the following conjecture: for any $0<\a<1$ and $0<\b<\a$, there is $C(\a,\b)>0$ such that for any $m$,
$$
\#\big\{d|m: m^{\a}\leq d\leq m^{\a}+m^{\b}\big\}\leq C(\a,\b).
$$   	
It turns out that this conjecture implies the mentioned conjecture (\ref{1.3}) of Cilleruelo and C\'ordoba (in fact, it implies Conjecture 1.5 from \cite{G} from which (\ref{1.3}) follows; see \cite{G} for details). It can also be shown by the same method that it implies the analog of (\ref{1.3}) for cubes as well, that is, it would improve exponent $2/3-\e$ in our main result to $1-\e$.
\end{rem}

\section{Proof of Theorem \ref{th1.2}}

For the first claim of the theorem, we note that if 
$$
(N+s_1)^3+(N+s_2)^3 = (N+s_3)^3+(N+s_4)^3
$$	
for some $0\leq s_i\leq (0.5N)^{1/2}$, then 
$$
3N^2(s_1+s_2-s_3-s_4)+3N(s_1^2+s_2^2-s_3^2-s_4^2)= s_3^3+s_4^3-s_1^3-s_2^3.
$$
It is easy to see that in the case $s_1+s_2\neq s_3+s_4$ the first summand from the left-hand side dominates the two other terms and we have a contradiction. Now let $s_1+s_2=s_3+s_4$. If we denote $u_i=N+s_i$, then $u_1+u_2=u_3+u_4$ and $u_1^3+u_2^3=u_3^3+u_4^3$. Since $u^3+v^3=(u+v)((u+v)^2-3uv)$, we get $u_1u_2=u_3u_4$ and thus $\{u_1,u_2\}=\{u_3,u_4\}$, as desired.

\medskip 

Now we turn to the second claim of Theorem \ref{th1.2}. Following Ramanujan's example
$$
1^3 + 12^3 = 9^3 + 10^3,
$$
we look for the relations of the form 
$$
u_1^3 + u_2^3 = u_3^3 + u_4^3 = U
$$
with
\begin{equation*}
u_1 + u_2 = v-6, \quad u_3 + u_4 = v, \quad 4U = (v-6)v(v+9).	
\end{equation*}
We claim that there are infinitely many such examples. First, we mention that there are infinitely many  positive integers $X,Y$ with
\begin{equation}\label{3.1}  
7X^2 + 114 = Y^2
\end{equation}
(note that in Ramanjuan's example we have $X=1$, $Y=11$). It follows from consideration
$$ 
{\sqrt 7} X + Y = (\sqrt 7 + 11) (3\sqrt{7} + 8)^k =(X'\sqrt7+Y')(X_0\sqrt7+Y_0)^k
$$
for positive integers $k$: here $(X',Y')=(1,11)$ is a solution of our generalized Pell's equation (\ref{3.1}) and $(X_0,Y_0)=(3,8)$ is the fundamental solution of the corresponding Pell's resolvent $7X^2-Y^2=1$. We refer the reader to \cite[Chapters 3 and 4]{QDE} for a detailed overview of Pell's equation.

Now let $X,Y$ be large enough positive integers such that (\ref{3.1}) holds. We define
$$
u_1 = \frac{X^2-Y}2+6, \quad u_2 = \frac{X^2+Y}2+6, 
$$
$$ 
u_3 = \frac{X^2-X}2+9, \quad u_3 = \frac{X^2+X}2+9,
$$
and $v=X^2+18$. Then because of (\ref{3.1}) we have
$$
4(u_1^3+u_2^3)=(X^2+12)(X^4+24X^2+3Y^2+144)=(X^2+12)(X^2+18)(X^2+27)=(v-6)v(v+9),
$$
and also
$$
4(u_3^3+u_4^3)=(v-6)v(v+9).
$$
Therefore,
$$
u_1^3 + u_2^3 = u_3^3 + u_4^3 
$$
and clearly $u_1,u_2,u_3,u_4$ are all positive integers of the from the interval $[N,N+CN^{1/2}]$, where $C>0$ is an absolute constant and 
$$
N=(X^2-Y)/2
$$ 
can be taken arbitrarily large since (\ref{3.1}) has infinitely many solutions. The claim follows.

\end{document}